\documentclass[12pt]{amsart}
\usepackage{diagram}
\usepackage{latexsym}
\input amssym.def
\input amssym.tex
\numberwithin{equation}{section}
\makeatletter
\renewcommand{\subsection}{\@startsection
{subsection}{2}{0mm}{\baselineskip}{-0.25cm}
{\normalfont\normalsize\rmfamily}}
\makeatother
\newtheorem{claim'}{Claim}
\newtheorem{theorem}{Theorem}[section]
\newtheorem{proposition}[theorem]{Proposition}
\newtheorem{lemma}[theorem]{Lemma}
\newtheorem{corollary}[theorem]{Corollary}
\newtheorem{scholium}[theorem]{Scholium}

\theoremstyle{remark}
\newtheorem*{remark}{Remark}
\newtheorem{example}[theorem]{Example}
\newenvironment{claim*}{\begin{trivlist}\item[\hskip%
\labelsep{\bf{Claim.}}]\it }%
{\end{trivlist}}

\def\d{{\mathcal D}}
\def\j{{\mathcal J}}
\def\pij{{\rm Fr}_{\mathcal J}}
\def\pix{{\rm Fr}_{X}}
\def\supp{{\rm Supp}}
\def\div{{\rm div}_\infty}

\def\kx{{k}(X)}
\def\fq{\mathbb F_{q^2}}
\begin{document}
\author[R.~Fuhrmann]{Rainer Fuhrmann}\thanks{The paper was partially written
while the first author was visiting the Instituto de Matem\'atica Pura e 
Aplicada, Rio de Janeiro (Oct. 1995 - Jan. 1996) supported by CNPq 
%
%
and the third author was visiting ICTP (Trieste - Italy) supported by the
International Atomic Energy Agency and UNESCO}
\author[A.~Garcia]{Arnaldo Garcia}
\author[F.~Torres]{Fernando Torres}
\title[Maximal curves]{On maximal curves}
\address{Fachbereich 6, Mathematik und Informatik, Universit\"at Essen, 
D-45117, Essen - Germany}
\email{RAINER.FUHRMANN@ZENTRALE.DEUTSCHE-BANK.dbp.de}
\address{IMPA, Estrada Dona Castorina 110, 22460-320, Rio de Janeiro -
Brazil}
\email{garcia@impa.br}
\address{Fachbereich 6, Mathematik und Informatik, Universit\"at Essen, 
D-45117, Essen - Germany}
\email{fernando.torres@uni-essen.de} 
\begin{abstract}
We study arithmetical and geometrical properties of {\it maximal curves}, 
that is, curves defined over the finite field $\mathbb F_{q^2}$ whose 
number of $\mathbb F_{q^2}$-rational points reaches the Hasse-Weil upper 
bound. Under a hypothesis on non-gaps at a rational point, we prove that 
maximal curves are $\mathbb F_{q^2}$-isomorphic to $y^q+y=x^m$, for some 
$m\in \mathbb Z^+$. As a consequence we show that a maximal curve of genus
$g=(q-1)^2/4$ is $\mathbb F_{q^2}$-isomorphic to the curve $y^q+y=x^{(q+1)/2}$.
\end{abstract}
\maketitle
\setcounter{section}{-1}
\section{Introduction} 
The interest on curves over finite fields was renewed after Goppa
\cite{Go} showed their 
applications to Coding Theory. One of the main features of linear codes
arising from curves is the fact that one can state a lower bound for their 
minimum distance. This lower bound is 
meaningful only if the curve has many rational points. The subject of this
paper is the study of {\it maximal curves}. 

Let $X$ be a projective, geometrically irreducible and 
non-singular algebraic curve 
defined over the finite field $\mathbb F_\ell$ with $\ell$ elements. A
celebrated 
theorem of Weil states that:
\[
\#\,X(\mathbb F_\ell) \le \ell+1+2g\sqrt\ell,
\]
where $X(\mathbb F_\ell)$ denotes the set of $\mathbb F_\ell$-rational points
of $X$ and $g$ is the genus of the curve. This bound was proved for
elliptic curves by Hasse. 

The curve $X$ is called \emph{maximal over} $\mathbb F_\ell$ (in this
case, $\ell$ must be a square; say $\ell = q^2$) if it attains the
Hasse-Weil upper bound; that is,
\[
\#\,X(\mathbb F_{q^2}) = q^2 + 1 + 2gq.
\]

Ihara \cite{Ih} shows that the genus of a maximal curve over $\mathbb
F_{q^2}$
satisfies:
\[
g \le (q-1)q/2.
\]

R\"uck and Stichtenoth \cite{R-Sti} show that the Hermitian curve (that
is, 
the curve given by $y^q+y = x^{q+1}$) is the unique (up to 
$\mathbb F_{q^2}$-isomorphisms) maximal curve over $\mathbb F_{q^2}$ having
genus $g = (q-1)q/2$.

It is also known that the genus of maximal curves over $\mathbb F_{q^2}$
satisfies (see \cite{F-T} and the remark after Theorem 1.4 here):
\[
g \le (q-1)^2/4 \qquad\text{or}\qquad g = (q-1)q/2\,.
\]

The Hermitian curve is a particular case of the following maximal
curves over $\mathbb F_{q^2}$\,:
\[
y^q + y = x^m, \text{ with $m$ being a divisor of $(q+1)$}.
\]

Note that the genus of the above curve is given by $g = (q-1)(m-1)/2$. 

%
%
%
In Section 1 we derive properties of 
maximal curves. The main tools being the application to the linear
system $\d = |(q+1)P_0|$, $P_0$ a rational point, of
St\"ohr-Voloch's approach \cite{S-V} to the Hasse-Weil bound via 
Weierstrass Point Theory; and Castelnuovo's genus bound for 
curves in projective spaces: \cite{C}, \cite[p. 116]{ACGH}, 
\cite[Corollary 2.8]{Ra}. A key result here is the fact that for any point
$P$ of the curve, the divisor $qP+\pix(P)$ is linearly equivalent to $\d$ 
(Corollary \ref{c1.2}). This is a
consequence of the particular fashion of the characteristic polynomial
$h(t)$ of the Frobenius endomorphism of the Jacobian of the curve, that
is, $h(t)$ is a power of a linear polynomial. This property also affects 
the geometry of the curve. More precisely, we show 
that maximal curves over $\mathbb F_{q^2}$ of genus $g \ge q-1$ are
non-classical curves for the  
canonical morphism (Proposition \ref{p1.7}). In some other cases one can 
deduce the non-classicality (for the canonical morphism) of the curve
from the knowledge of $h(t)$. We will see this for the Deligne-Lusztig
curve associated to the Suzuki group and to the Ree group (Proposition
\ref{p?}). The non-classicality of the curve corresponding to the Suzuki
group was already proved in \cite{G-Sti}. Our proof is different. It seems
that the curve corresponding to the Ree group provides a new example of a
non-classical curve. 

In Section 2, we characterize the curves
\[
y^q + y = x^m,\,\,m \text{ being a divisor of $(q+1)$},
\]
among the maximal curves over $\mathbb F_{q^2}$\,. This characterization being
in terms of non-gaps at a rational point (Theorem \ref{t2.3}). Finally in 
Section 3, applying the results of Section 2, we show that
\[
y^q + y = x^{(q+1)/2}\,, \text{ with $q$ odd},
\]
is the unique (up to $\mathbb F_{q^2}$-isomorphisms) maximal curve over
$\mathbb F_{q^2}$ with $g = (q-1)^2/4$. 
\section{Maximal curves}
Throughout this paper we use the following notation: 
\begin{itemize}
\item By a curve we mean a projective, geometrically irreducible,   
non-singular algebraic curve defined over a finite field. 
\item Let $k$ denote the finite field with $q^2$ elements, where 
$q$ is a power of a prime $p$. Let ${\bar k}$ denote its algebraic closure.
\item The symbol $X(k)$ (resp. $\kx$) stands for the set 
of $k$-rational points (resp. for the field of $k$-rational functions) of 
a curve $X$ defined over $k$.
\item If $x\in \kx$, then ${\rm div}(x)$ (resp. $\div(x)$) denotes 
the divisor (resp. the pole divisor) of the function $x$.
\item Let $P$ be a point of a curve. Then $v_P$ (resp. $H(P)$) stands 
for the valuation (resp. for the Weierstrass non-gap semigroup) associated to
$P$. We  denote by $m_i(P)$ the $i$th non-gap at $P$.
 \item Let $D$ be a divisor on $X$ and $P\in X$. We denote by 
${\rm deg}(D)$ the degree of $D$, by $\supp(D)$ the support of $D$, and by 
$v_P(D)$ the coefficient of $P$ in $D$. If $D$ is a $k$-divisor, we set
$$
L(D):= \{f\in \kx \mid {\rm div}(f)+D \succeq 0\}, \quad\text{and}\quad
\ell(D):={\rm dim}_k\, L(D).
$$
\item The symbol ``$\sim$" denotes linear equivalence of divisors.
\item The symbol $g^r_d$ stands for a linear system of projective dimension $r$
and degree $d$. 
\end{itemize}

We first review some facts from Weierstrass Point Theory (see \cite{Sch}  
and \cite{S-V}).

\noindent \textbf{Weierstrass points.}\label{wp} 
Let $X$ be a curve of genus $g$, and $\d=g^r_d$ be a base-point-free $k$-linear
system on $X$. Then associated to a point $P\in X$ we have the Hermitian 
$P$-invariants $j_0(P)=0<j_1(P)<\ldots<j_r(P)\le d$ of $\d$ (also called 
the $(\d,P)$-orders). This sequence is the same for all but finitely many 
points. These finitely many points $P$, where exceptional $(\d,P)$-orders 
occur, are called the $\d$-Weierstrass points of $X$. The Weierstrass
points of the curve are those exceptional points obtained from the canonical
linear system. A curve is called {\it non-classical} if the generic order
sequence (for the canonical linear system) is different from
$\{0,1,\ldots,g-1\}$.

Associated to the linear system $\d$ there exists a divisor $R$ supporting
exactly the $\d$-Weierstrass points. Let 
$\epsilon_0<\epsilon_1<\ldots<\epsilon_r$ denote the $(\d,Q)$-orders for 
a generic point $Q\in X$. Then we have
\begin{equation}\label{ineq1}
\epsilon_i \le j_i(P),\quad \text{for each}\quad i=0,1,2,\dots,r \text{ and for
any point $P$,}
\end{equation}
and also that
\begin{equation}\label{degR}
{\rm deg}(R)= (\epsilon_1+\ldots+\epsilon_r)(2g-2)+(r+1)d.
\end{equation}
Associated to $\d$ we also have a divisor $S$ whose support contains the set
$X(k)$ of $k$-rational points on $X$. Its degree is given by 
\begin{equation*}
{\rm deg}(S)=(\nu_1+\ldots+\nu_{r-1})(2g-2)+(q^2+r)d,
\end{equation*}
where the $\nu_i's$ form a subsequence of the $\epsilon_i's$. More 
precisely, there exists an integer $I$ with $0<I\le r$ such that 
$\nu_i=\epsilon_i$ for $i<I$, and $\nu_i=\epsilon_{i+1}$ otherwise. 
Moreover, for $P\in X(k)$, 
\begin{equation}
v_P(S)\ge \sum_{i=1}^{r}(j_i(P)-\nu_{i-1}),\ \text{and}\quad \nu_i\le
j_{i+1}(P)-j_1(P),\ \text{for each $i=1,2,\dots,r$}.
\end{equation}

\noindent \textbf{Maximal curves.} We study some arithmetical and
geometrical properties of maximal curves. To begin with we recall the
following basic result concerning Jacobians. Let $X$ be a curve, $\pij$
the Frobenius endomorphism (relative to the base field) of the Jacobian
$\j$ of $X$, and $h(t)$ the characteristic polynomial of $\pij$. Let $
h(t)=\prod_{i=1}^{T}h^{r_i}_i(t)$ 
be the factorization over $\mathbb Z[t]$ of $h(t)$. Then 
\begin{equation}\label{car}
\prod_{i=1}^{T}h_i(\pij)=0\qquad \mbox{on}\ \ \j.
\end{equation}
This follows from the semisimplicity of $\pij$ and the fact that the
representation of endomorphisms of $\pij$ on the Tate module is faithful
(cf. \cite[Thm. 2]{Ta}, \cite[VI, \S3]{L}).

In the case of a maximal curve over $k = \mathbb F_{q^2}$,
$h(t)=(t+q)^{2g}$. Therefore from (\ref{car}) we obtain the following 
result, which is contained in the proof of \cite[Lemma 1]{R-Sti}.
\begin{lemma}\label{l1.1} 
The Frobenius map $\pij$ 
(relative to $k$) of the Jacobian $\j$ of a maximal curve over $k$ acts as
multiplication by
$(-q)$ on $\j$.
\end{lemma}

Let $X$ be a maximal curve over $k$. Fix $P_0\in X(k)$, and consider the
map 
$f=f^{P_0}: X\to \j$ given by $P\to [P-P_0]$. We have
\begin{equation}\label{car1}
f\circ {\rm Fr} = \pij\circ f,
\end{equation}
where Fr denotes the Frobenius morphism on $X$ relative to $k$. Hence,
from (\ref{car1}) and Lemma \ref{l1.1}, we get:
\begin{corollary}\label{c1.2}
For a maximal curve $X$ over $k$, it holds
$$
{\rm Fr }(P)+qP \sim (q+1)P_0, \,\,\text{ for all points $P$ on $X$.}
$$
\end{corollary}

It follows then immediately that
\begin{corollary}(\cite[Lemma 1]{R-Sti})\label{c1.3}   
Let $X$ be a maximal curve over $k$, $P_0, P_1 \in X(k)$. Then
$(q+1)P_1\sim (q+1)P_0$.
\end{corollary}

Consider now the linear system $\d = g^{n+1}_{q+1}:= 
|(q+1)P_0|$. Corollary \ref{c1.3} says that $\d$ is a $k$-invariant of 
the curve. In particular, its dimension $n+1$ is independent of the choice of
$P_0\in X(k)$. Moreover from Corollary \ref{c1.3} we have that $q+1\in
H(P_0)$; i.e., $(q+1)$ is a non-gap at a rational point, and hence $\d$ is
base-point-free. From now on the letter $\d$ will always denote the linear
system $|(q+1)P_0|$, \,\,$P_0$ a rational point, $(n+1)$ being its projective
dimension, $R$ will always mean the divisor supporting exactly the 
$\d$-Weierstrass points, and Fr will always stand for the Frobenius
morphism on $X$ relative to $k$.
\begin{theorem}\label{t1.4} For a maximal curve $X$ over $k$, the $\d$-orders
satisfy (notations being as above):
\begin{enumerate}
\item[(i)] $\epsilon_{n+1}=\nu_n=q$.
\item[(ii)] $j_{n+1}(P)=q+1$ if $P\in X(k)$, and $j_{n+1}(P)=q$ otherwise; in
particular, all rational points over $k$ are $\d$-Weierstrass points of $X$.
\item[(iii)] $j_1(P)=1$ for all points $P\in X$; in particular,
$\epsilon_1=1$.
\item[(iv)] If $n\ge2$, then $\nu_1 = \epsilon_1 = 1$.
\end{enumerate}
\end{theorem}
\begin{proof} Statement (iii), for $P\in X(k)$,
follows from (i), (ii) and 
the second inequality in (1.3). From Corollary \ref{c1.2} it follows the
assertion
(ii) and $\epsilon_{n+1}=q$. Furthermore, it also follows that $j_1(P)=1$ for 
$P\not\in X(k)$: in fact, let $P'\in X$ be such that ${\rm Fr }(P')=P$; then
$P+qP'= {\rm Fr }(P')+ qP'\sim (q+1)P_0$.

Now we are going to prove that $\nu_n=\epsilon_{n+1}$. Let $P\in 
X\setminus \{P_0\}$. Corollary \ref{c1.2} says that $\pi({\rm Fr }(P))$ belongs
to the osculating hyperplane at $P$, where $\pi$ stands for the morphism
associated to $\d$. This morphism $\pi$ can be defined by a base
$\{f_0,f_1,\ldots,f_{n+1}\}$ of\linebreak $L((q+1)P_0)$. Let $x$ be a
separating variable of $k(X)\mid k$. Then by \cite[Prop. 1.4(c),  
Corollary 1.3]{S-V} the rational function below is identically zero
$$
w:= {\rm det} \begin{pmatrix}
f_0\circ{\rm Fr}  & \ldots  &  f_{n+1}\circ{\rm Fr}\\
D^{\epsilon_0}_x f_0 & \ldots & D^{\epsilon_0}_x f_{n+1}\\
\vdots      & \vdots &  \vdots\\
D^{\epsilon_n}_x f_0 & \ldots & D^{\epsilon_n}_x f_{n+1}
\end{pmatrix}\,\,,
$$
since it satisfies $w(P)=0$ for a generic point $P$. Let $I$ be the
smallest integer such that the row $(f_0\circ{\rm Fr},\ldots,f_{n+1}\circ{\rm
Fr})$ is a linear combination of the vectors $(D^{\epsilon_i}_x 
f_0,\ldots,D^{\epsilon_i}_x f_{n+1})$ with $i=0,\ldots,I$. Then according 
to \cite[Prop. 2.1]{S-V} we have
$$
\{\nu_0<\ldots<\nu_n\}=\{\epsilon_0<\ldots<\epsilon_{I-1}
<\epsilon_{I+1}<\ldots<\epsilon_{n+1}\}.
$$

That $\epsilon_1=1$ follows from statement (iii). Suppose that $\nu_1 > 1$.
Since $j_1(P)=1$ for all points of $X$, it follows from the proof of
\cite[Thm. 1]{H-V} that
\[
\#\,X(k) = (q+1)(q^2-1) - (2g-2).
\]
From the maximality of $X$, we then conclude $2g = (q-1)\cdot q$.

\noindent On the other hand, $\pi$ is a birational morphism as 
follows from \cite[Prop. 1]{Sti-X} (see also Proposition \ref{p1.5}(iv) 
here). Then Castelnuovo's genus bound for curves in projective spaces 
applied to  
the morphism $\pi$ reads: 
\begin{equation}\label{eq1.4}
2g \le M\cdot(q-n+e) \le
\begin{cases}
(2q-n)^2/4n\,, &\quad \text{ if $n$ is even}\\
((2q-n)^2-1)/4n\,, &\quad \text{ if $n$ is odd,}
\end{cases}
\end{equation}
where $M$ is the integer part of $q/n$ and $e = q-M\cdot n$. We then
conclude that $n=1$ and this finishes the proof of the
theorem.
\end{proof} 
\begin{remark}
For a maximal curve $X$ with $n=1$, we have $\nu_1 = \epsilon_2 = q > 1$. Then
the proof above shows that $2g = (q-1)\cdot q$. It then follows from
\cite{R-Sti} that the curve $X$ is $k$-isomorphic to the Hermitian curve
given by $y^q+y = x^{q+1}$. Also, if $n \ge 2$ then from Castelnuovo's
formula (\ref{eq1.4}) we get $g \le (q-1)^2/4$. This is the main result of
\cite{F-T}.
\end{remark}
The next proposition gives information on $\d$-orders and non-gaps at
points of $X$.
\begin{proposition}\label{p1.5}
Let $X$ be a maximal curve over $k$ (notations being as before). Then:
\begin{enumerate}
\item[(i)] For each point $P$ on $X$, we have $\ell(qP) = n+1$; i.e., we have
the following behaviour for the non-gaps at $P$
\[
0 < m_1(P) <\cdots< m_n(P) \le q < m_{n+1}(P).
\]
\item[(ii)] If $P$ is not rational over $k$, the numbers
below are $\d$-orders at the point $P$
\[
0 \le q-m_n(P) < \cdots < q-m_1(P) < q.
\]
\item[(iii)] If $P$ is rational over $k$, the numbers below are exactly the
$(\d,P)$-orders 
\[
0 < q+1-m_n(P) <\cdots< q+1 - m_1(P) < q+1.
\]
In particular, if $j$ is a $\d$-order at a rational point $P$ then $q+1-j$ is a
non-gap at $P$. 
\item[(iv)] If $P \in X(\mathbb F_{q^4})\backslash X(k)$, then $q-1$ is a 
non-gap at $P$. If $P\notin X(\mathbb F_{q^4})$, then $q$ is a non-gap at $P$.
If $P \in X(k)$, then $q$ and $q+1$ are non-gaps at $P$. 
\item[(v)] Let $P$ be a non-Weierstrass point of $X$ (for the canonical
morphism) and suppose that $n\ge2$, then we have for the non-gaps at $P$ that
$m_{n-1}(P) = q-1$ and $m_n(P)=q$.
\end{enumerate}
\end{proposition}
\begin{proof}
Assertion (i) follows from Corollary 1.2. Let $m(P)$ be a non-gap at a point
$P$ of $X$ with $m(P) \le q$, then by definition there exists a positive
divisor $E$ disjoint from $P$ with
\[
E \sim m(P)\cdot P.
\]
Summing up to both sides of the equivalence above the divisor $(q-m(P))\cdot
P+{\rm Fr }(P)$, we get
\[
E + (q-m(P))\cdot P + {\rm Fr }(P) \sim qP + {\rm Fr }(P) \sim (q+1)P_0\,.
\]
This proves assertions (ii) and (iii). To prove assertion (iv) we just 
apply (as in \cite[IV, Ex.
2.6]{Har}) the Frobenius morphism to the
equivalence in Corollary \ref{c1.2}, getting
\[
{\rm Fr}^2(P) + (q-1){\rm Fr }P \sim qP.
\]
The fact that $q$ and $q+1$ are non-gaps at any rational point follows from
assertion (iii) taking $j=0$ and $j=1$.

\noindent Now we are going to prove the last assertion (v). From assertion (iv)
we know already 
\[
m_n(P) = q \quad{\rm and}\quad m_{n-1}(P) \le q-1.
\]

Suppose that $m_{n-1}(P) < (q-1)$. It then follows from Theorem 1.4 and the
assertion (ii) above that the generic order sequence for the linear system $\d$
is as given below:
\[
\epsilon_0=0 < \epsilon_1=1 < \epsilon_2 = q-m_{n-1}(P) <\cdots< \epsilon_n =
q-m_1(P) < \epsilon_{n+1} = q.
\]
On the other hand, we have that Equation (1.1) implies
\[
m_i(Q) \le m_i(P), \text{ for each $i$ and each $Q\in X$}.
\]
Thus at a rational point $Q \in X$, it follows from assertion (iii) that:
\[
v_Q(R) \ge \sum_{i=1}^{n+1}(j_i(Q)-\epsilon_i) = 1 + \sum_{i=1}^{n-1}(m_i(P)-
m_i(Q)+1) \ge n.
\]
From the maximality of $X$, Equation (1.2) and \cite[Thm. 1]{Ho}, we
conclude that
\[
n(q^2+2gq+1) \le {\rm deg }R \le (n+2)\epsilon_{n+1}(g-1) + (n+2)(q+1).
\]
Using that $\epsilon_{n+1} = q$, we finally have $nq^2 + qg(n-2) \le 2$. 
This contradicts the assumption that $n \ge 2$.
\end{proof}

\begin{example}\label{e1.6}
By Theorem \ref{t1.4} we have that all rational points of the curve are
$\d$-Weierstrass points. However, these sets may be different from each
other as the following example shows:

Let $X$ be the hyperelliptic curve defined by $x^2+y^5=1$ over $\mathbb
F_{81}$. The curve $X$ is maximal because it is covered by the Hermitian curve
$x^{10}+y^{10}=1$  (see \cite[Example VI.4.3]{Sti}). It has genus 2 and at 
a
generic point $P$, we have 
$m_7(P)=9$. Hence we have $\d=|10P_0|=g^8_{10}$. All the canonical Weierstrass
points are trivially rational points, and since $\#X(\mathbb F_{81})=118 >
\#\,\{\text{Weierstrass points}\} = 6$, we have two possibilities for the
$(\d,P)$-orders at rational points, namely:

(a) If $P$ is a rational non-Weierstrass point; then its orders are
$0,1,2,3,4,5,6,7,10$.

(b) If $P$ is a Weierstrass point; then its orders are 
$0,1,2,3,4,5,6,8,10$.

These computations follow from Proposition 1.5(iii). From the $\d$-orders in
(a) above, we conclude that the generic order sequence for $\d$ is
$0,1,2,3,4,5,6,7,9$. Hence, ${\rm deg}(R)=164$ and $v_P(R)=1$ (resp.
$v_P(R)=2$) if $P$ satisfies (a) (resp. (b)) above. Since ${\rm
deg }R -112\times1-6\times2=40>0$, we then  conclude that there exist
non-rational $\d$-Weierstrass points. The order sequence at such points must
be $0,1,2,3,4,5,6,8,9$ and so there exist 40 non-rational $\d$-Weierstrass
points, namely the fixed points of $\sigma\circ{\rm Fr}$, where $\sigma$
denotes the hyperelliptic involution.
\end{example}
By Proposition \ref{p1.5}(v) we have that $q-n$ is a lower  
bound for the genus of a maximal curve over $\fq$. We are going to show
that classical (for the canonical morphism) maximal curves attain such a
bound.  
\begin{proposition}\label{p1.7}
Let $X$ be a maximal curve over $k = \mathbb F_{q^2}$ and let $g \ge 2$ be its
genus. Then
\begin{enumerate}
\item[(i)] If $g > q-n$ (with $n+1=\dim \d$ as before), then $X$ is 
non-classical for the canonical morphism. In particular, this holds for $g
\ge q-1$.
\item[(ii)] If $X$ is hyperelliptic and the characteristic is two, then $X$ has
just one Weierstrass point for the canonical morphism.
\end{enumerate}
 \end{proposition}
\begin{proof}
(i) If $X$ is classical, then at a generic point $P$ of the curve $X$ we have
\[
m_i(P) = g+i, \quad \forall\,i \in \mathbb N.
\]
On the other hand, from Proposition 1.5(iv) we have
$m_n(P) = q$. We then conclude that $g+n=q$. Now if $g\ge q-1$ and $X$ 
classical, then $n=1$. Therefore from the remark after Theorem 1.4 we
would have $2g=q(q-1)$ and so $g=1$, a contradiction.

(ii) Since $X$ is hyperelliptic, the Weierstrass points are the fixed
points of the hyperelliptic involution. Let $P, Q$ be Weierstrass points
of $X$ (they exist because the genus is bigger than one). From
$2P\sim 2Q$ and $2\mid q$, we get $qP\sim qQ$. Therefore by Corollary 
\ref{c1.2},
$$
qQ+\text{ Fr}(Q)\sim qP+\text{ Fr}(P)\sim qQ+\text{ Fr}(P),
$$
and so ${\rm Fr}(P)\sim {\rm Fr}(Q)$. This implies ${\rm Fr}(P)={\rm Fr}(Q)$,
and consequently $P=Q$.
\end{proof}
%
%
%
%
%
\begin{remark}
Hyperelliptic maximal curves are examples of classical curves for the 
canonical morphism. It would be interesting to investigate
the maximal curves that are both non-hyperelliptic
and classical for the canonical morphism. Examples of such curves are the 
one of genus 3 over $\mathbb F_{25}$ listed by Serre in 
\cite[\S4]{Se}, and the generalizations of Serre's example obtained by
Ibukiyama \cite[Thm. 1]{I}.

Another question is whether or not the condition $g=q-n$ characterizes
classical (for the canonical morphism) maximal curves.
\end{remark}

Now we present two non-classical (for the canonical morphism) maximal
curves over $\fq$ of genus 
$g<q-1$. These are the so-called Deligne-Lusztig curves associated to the
Suzuki group and to the Ree group.
\begin{proposition}\label{p?} (I) Let $s\in \mathbb N,\ r:= 2^{2s+1},\ 
r_0:= 2^s$, and consider the curve $X$ over $\mathbb F_r$ defined by
$$
y^r-y=x^{r_0}(x^r-x).
$$
Then
\begin{enumerate}
\item[(i)] (\cite{H}, \cite{H-Sti}, \cite{Se}) The genus of $X$ is
$g=r_0(r-1)$ and this curve is maximal over $\mathbb F_{r^4}$.
\item[(ii)] (\cite{G-Sti}) The curve $X$ is non-classical for the
canonical morphism.
\end{enumerate}
(II) Let $s\in N,\ r:=3^{2s+1},\ r_0:=3^s$, and consider the curve $X$
over
$\mathbb F_r$ defined by
$$
y^r-y=x^{r_0}(x^r-x),\qquad z^r-z=x^{2r_0}(x^r-x).
$$
Then
\begin{enumerate}
\item[(i)] (\cite{H}, \cite{P}, \cite{Se}) The genus of $X$ is
$g=3r_0(r-1)(r+r_0+1)$ and this curve is maximal over $\mathbb F_{r^6}$.
\item[(ii)] The curve $X$ is non-classical for the canonical morphism.
\end{enumerate}
\end{proposition}
\begin{proof} We first set some notations. We write ${\rm Fr}$ for the
Frobenius morphism on $X$ relative to $\mathbb F_{r}$ and $h_i(t)$ for
the characteristic polynomial of the Frobenius endomorphism (relative to
$\mathbb F_{r^i}$) of the Jacobian of $X$.

(I) From \cite[Prop. 4.3]{H}, \cite{H-Sti}, \cite{Se} we
know that $g=r_0(r-1)$ and that $h_1(t)=(t^2+2r_0t+r)^g$. If $a_1$ and
$a_2$ denote the roots of $h_1(t)$, then we have $a_1+a_2=-2r_0$ and
$a_1a_2=r$. It then follows easily that $(a_1a_2)^4=r^4$ and
$a_1^4+a_2^4=-2r^2$, and hence that $h_4(t)=(t+r^2)^{2g}$. This shows the
maximality over $\mathbb F_{r^4}$. 

\noindent Now by (\ref{car}) we have 
$\pij^2+2r_0\pij+r=0$ on $\j$ and then 
by (\ref{car1}) we obtain 
$$
{\rm Fr}^2(P)+2r_0{\rm Fr}(P)+rP\sim (1+2r_0+r)P_0,
$$
for all $P\in X$, $P_0\in X(\mathbb F_r)$. Now applying ${\rm Fr}$ to the
equivalence above we get 
$$
{\rm Fr}^3(P)+(2r_0-1){\rm Fr}^2(P)+(r-2r_0){\rm Fr}(P)\sim rP.
$$
Hence we conclude that $r\in H(P)$ at a general point $P\in X$, and since
$g\ge r$, that $X$ is non-classical for the canonical morphism.

(II) From \cite[Prop. 5.3]{H}, \cite{P}, \cite{Se} we already know the
formula for 
the genus and that $h_1(t)=(t^2+r)^a(t^2+3r_0t+r)^b$ with $a, 
b\in \mathbb N$ and $a+b=g$. Let $a_1, a_2$ denote the roots of some 
factor of $h_1(t)$. Then in either case we get $(a_1a_2)^6=r^6$ and
$a_1^6+a_2^6=-2r^3$ and hence that $h_6(t)=(t+r^3)^{2g}$. This shows the
maximality of $X$ over $\mathbb F_{r^6}$. 

Finally as in the proof of item (I) we conclude
that $r^2\in H(P)$ for a generic $P\in X$. Since $g\ge r^2$ the
assertion follows.
\end{proof}

To finish this section on maximal curves, we study some properties 
involving the morphism $\pi: X\to \mathbb P^{n+1}$ associated to the linear
system $\d=|(q+1)P_0|$. 

\begin{proposition}\label{p1.8}
The following statements are equivalent:
\begin{enumerate} 
\item[(i)] The morphism $\pi$ is a closed embedding, 
i.e. $X$ is $k$-isomorphic to $\pi(X)$.
\item[(ii)] For all $P \in X(\mathbb F_{q^4})$, we have that $\pi(P)\in \mathbb
P^{n+1}(k) \Leftrightarrow P \in X(k)$.
\item[(iii)] For all $P \in X(\mathbb F_{q^4})$, we have that $q$ is a non-gap
at $P$.
\end{enumerate}
\end{proposition}
\begin{proof}
Let $P\in X$. Since $j_1(P)=1$ (cf. Theorem \ref{t1.4}(iii)), we have 
that $\pi(X)$ is non-singular at all branches centered at 
$\pi(P)$. Thus $\pi$ is an embedding if and only if $\pi$ is injective.
\begin{claim*}
We have $\pi^{-1}(\pi(P))\subseteq \{P,{\rm Fr}(P)\}$ and if $\pi$ is not
injective at the point $P$, then $P\in X(\mathbb F_{q^4})\setminus 
X(k)$ and $\pi(P)\in \mathbb P^{n+1}(k)$.
\end{claim*}

From Corollary \ref{c1.2} it follows that $\pi^{-1}(\pi(P))\subseteq 
\{P,{\rm Fr}(P)\}$. Now if $\pi$ is not injective at $P$, then $P\not\in
X(k)$ and, since $P\in \pi^{-1}(\pi({\rm Fr}(P)))\subseteq 
\{{\rm Fr}(P), {\rm Fr}^2(P)\}$, we have ${\rm Fr}^2(P)=P$, i.e. $P\in
X(\mathbb F_{q^4})\setminus X(k)$. Furthermore we have 
$\pi(P)=\pi({\rm Fr}(P))={\rm Fr}(\pi(P))$, i.e. $\pi(P)\in \mathbb
P^{n+1}(k)$. This proves the claim above. 

From this claim the equivalence (i) $\Leftrightarrow$ (ii) follows 
immediately. As to the implication (i) $\Rightarrow$ (iii),  
we know that $\dim |{\rm Fr}(P)+qP-P-{\rm Fr}(P)|=\dim |{\rm Fr}(P)+qP| - 2$ 
(Corollary \ref{c1.2} and \cite[Prop. 3.1(b)]{Har}), i.e. $\ell((q-1)P)=n$
and so $q\in H(P)$, by Proposition 1.5(i).
 
Finally we want to conclude from (iii) that $\pi$ is an embedding.  
According to the above claim it is sufficient to show that 
$\pi^{-1}(\pi(P))=\{P\}$, for $P \in X(\mathbb F_{q^4})\setminus X(k)$. Let
then $P\in X(\mathbb F_{q^4})$. Since we have $q\in H(P)$, there is a    
divisor $D\in |qP|$ with $P\notin \supp(D)$. In particular,
$$
{\rm Fr}(P)+D\sim {\rm Fr}(P)+qP\sim (q+1)P_0,
$$
and then $\pi^{-1}(\pi({\rm Fr}(P)))\subseteq \supp({\rm Fr}(P)+D)$. So if
$\pi^{-1}(\pi(P))=\{P,{\rm Fr}(P)\}$, then we would have that $P\in \supp(D)$,
a contradiction. This means altogether that $\pi$ is injective 
and so indeed a closed embedding.
\end{proof}
\begin{remark}
Condition (iii) above is satisfied whenever $q\ge 2g$, and in most of the
well known examples of maximal curves the morphism $\pi$ is always an
embedding. Then a natural question is whether or not $\pi$ is an embedding
for an arbitrary maximal curve. We conjecture that this property is a 
necessary condition for a maximal curve being covered by the Hermitian curve. 
\end{remark}
\begin{proposition}\label{p1.9} 
Suppose that $\pi: X\rightarrow \mathbb P^{n+1}$ is a closed embedding.
Let $P_0\in X(k)$ and assume furthermore that there exist $r, s \in
H(P_0)$ such that all non-gaps at $P_0$ less than or equal to $q+1$ are
generated by $r$ and $s$. Then the semigroup $H(P_0)$ is generated by $r$ and
$s$. In particular, the genus of $X$ is equal to $(r-1)(s-1)/2$.
\end{proposition}
\begin{proof}
Let $x, y\in k(X)$ with $\div(x)=sP_0$ and $\div(y)=rP_0$. Since we have that
$q, q+1 \in H(P_0)$, then the numbers $r$ and $s$ are coprime. Let  $\pi_2: 
X\rightarrow \mathbb P^2$, be given by $P\mapsto (1:x(P):y(P))$. Then the
curves $X$ and $\pi_2(X)$ are birational and the image $\pi_2(X)$ 
is a plane curve given by an equation of the type below:
$$
x^r+\beta y^s+\sum_{is+jr<rs} \alpha_{ij}x^iy^j=0,
$$
where $\beta,\alpha_{ij}\in k$ and $\beta\neq 0$. We are going to prove that 
$\pi_2(P)$ is a non-singular point of the curve $\pi_2(X)$ for all $P\neq P_0$.
From this it follows that $g=(r-1)(s-1)/2$ and also that $H(P_0)=\langle 
r,s\rangle$ (see \cite[Ch. 7]{Ful}, \cite{To}). 
 
Let $1,f_1,\ldots,f_{n+1}$ be a basis of $L((q+1)P_0)$, where 
$n+1:=\dim |(q+1)P_0|$. Then there exist polynomials $F_i(T_1,T_2)\in 
k[T_1,T_2]$ for $i=1,\ldots,n+1$, such that  
$$ 
f_i=F_i(x,y), \qquad \mbox{for}\qquad i=1,\ldots, n+1.
$$
The existence of these polynomials follows from the hypothesis on the non-gaps
at $P_0$ less than or equal to $(q+1)$.

\noindent Consider the maps $\pi |(X\setminus\{P_0\}):
X\setminus\{P_0\}\rightarrow \mathbb  A^{n+1}$ given by $P\mapsto
(f_1(P),\ldots,f_{n+1}(P))$; 
$\pi_2 |(X\setminus\{P_0\}): X\setminus\{P_0\}\rightarrow \mathbb A^2$,  
$P\mapsto (x(P),y(P))$; and $\phi :  \mathbb A^2\rightarrow \mathbb 
A^{n+1}$, given by $(p_1,p_2)\mapsto 
(F_1(p_1,p_2),\ldots,F_{n+1}(p_1,p_2))$.
Then the following diagram is commutative
$$
\Atriangle[X\setminus\{P_0\}`{\mathbb A}^2`{\mathbb A}^r;\pi_2`\pi`\phi].
$$
Thus we have for a point $P$ of
$X\setminus\{P_0\}$ and the corresponding local rings assigned to 
$\pi(P), \pi_2(P)$ the commutative diagram      
$$
\Atriangle[O_{\pi(X),\pi(P)}`O_{\pi_2(X),\pi_2(P)}`O_{X,P};f`c`h]\ ,
$$
where $h$ is injective since $k(X)=k(x,y)$, and $c$ is an isomorphism by 
assumption. Thus $\pi_2{X}$ is non-singular at $\pi_2{P}$.
\end{proof}
\section{Certain maximal curves}
The curves we have in mind in this section are the ones given by (see 
\cite{G-V} 
and \cite{Sch}):
\[
y^q + y = x^m,\,\, \text{ with $m$ being a divisor of $(q+1)$}.
\]
These are maximal curves (with $2g=(m-1)(q-1)$) since they are covered by the
Hermitian curve. If $P_0$ is the unique point at 
infinity of this curve, then the semigroup of non-gaps at $P_0$ is generated by
$m$ and $q$ and we have:
\[
m\cdot n = q+1, \quad\text{where}\quad (n+2) = \ell((q+1)P_0). \qquad (**)
\]
The goal of this section is to give a proof that the above condition
$(**)$ on
non-gaps at a rational point $P_0$ characterizes the curves $y^q+y = x^m$ among
the maximal curves over the finite field $k$.

\begin{proposition}\label{p2.1} Let $X$ be a maximal curve of genus $g$.  
Suppose that there exists a rational point $P_0\in X(k)$ such that $n\cdot
m=q+1$, with $m$ being a non-gap at $P_0$\,. Then, we have $2g=(q-1)(m-1)$.
Also, there are at most two types of $(\d,P)$-orders at rational points $P \in
X(k)$:

{\bf Type 1.} The $\d$-orders at $P$ are\,\,
$0,1,2,3,\dots,n,q+1$. In this case we have $v_P(R)=1$.

{\bf Type 2.} The $\d$-orders at $P$ are\,\,
$0,1,m,2m,\dots,(n-1)m,
q+1$. In this case we have $w_2 := v_P(R) = n((n-1)m-n-1)/2\,\, + 2$.

\noindent Moreover, the set of $\d$-Weierstrass points of $X$ coincides
with the set of its $k$-rational points, and the order sequence for $\d$
is $0,1,2,\dots,n,q$.
\end{proposition}
\begin{proof} The morphism $\pi$ can be 
defined by $(1:y:\ldots:y^{n-1}:x:y^n)$, where $x, y \in k(X)$ are functions
such that 
$$
\div(x)=qP_0\qquad {\rm and}\qquad \div(y)=mP_0.
$$
The set of $\d$-orders at $P_0$ is of Type 2, as follows from 
Proposition 1.5(iii).

\noindent Let $P\in X\setminus \{P_0\}$. From the proof of \cite[Thm. 
1.1]{S-V} and letting $z=y-y(P)$, we have 
\begin{equation}\label{vals}
v_P(z),\ldots,nv_P(z)
\end{equation} 
are $(\d,P)$-orders. Thus, considering a non-ramified point for $y:X\to \mathbb
P^1$, we conclude that the order sequence of the linear system $\d$ is given by
$$
\epsilon_i = i\qquad {\rm for}\ \ i=1,\ldots, n, \quad {\rm and}\ \
\epsilon_{n+1}=q.
$$

If $P$ is a rational point, by Theorem 1.4, we know that 1 and $(q+1)$ are
$(\d,P)$-orders. We consider two cases:
\begin{enumerate}
\item[(1)] $v_P(z) = 1$: This implies that the point $P$ is of Type 1.
\item[(2)] $v_P(z) > 1$: From assertion (2.2) above, it follows $n\cdot v_P(z)
= q+1$ and hence $v_P(z)=m$. Then, we have that the point $P$ is of Type 2.
\end{enumerate}

If $P$ is not a rational point, by Theorem 1.4, we have that $j_{n+1}(P)=q$. If
$v_P(z) > 1$ and using assertion (\ref{vals}), we get
\[
n\cdot v_P(z) = q = n\cdot m-1.
\]
Hence $n=1$ and the $(\d,P)$-orders are $0,1,q$. This shows that $P$ is not a
$\d$-Weierstrass point. If $v_P(z)=1$, again from assertion (\ref{vals}), 
we have that
\[
0,1,2,\ldots,n,\,\,q
\]
are the $\d$-orders at the point $P$; i.e., $P$ is not a $\d$-Weierstrass
point. This shows the equality of the two sets below:
\[
\{\d-\text{Weierstrass points of } X\} = \{k-\text{rational points of } X\}.
\]
The assertions on $v_P(R)$ follow from \cite[Thm. 1.5]{S-V}. 

Let $T_1$ (resp. $T_2$)  denote the number of rational points $P\in X(k)$ whose
$(\d,P)$-orders are of Type 1 (resp. Type 2). Thus we have from the equality in
(1.2)
$$
{\rm deg}(R)= (n(n+1)/2+q)(2g-2)+(n+2)(q+1)= T_1 + w_2 T_2.
$$
Riemann-Hurwitz applied to $y:X\to \mathbb P^1({\bar k})$ gives  
$$
2g-2=-2m + (m-1)T_2\,.
$$
Since $T_1+T_2 = \#X(k)= q^2+2gq+1$, and using the two equations above, we
conclude after tedious computations that $2g=(m-1)(q-1)$. This finishes the
proof of the proposition.
\end{proof}

Now we are going to prove that maximal curves as in Proposition 2.1 are
isomorphic to $y^q+y = x^m$. To begin with we first generalize \cite[Lemma 
5]{R-Sti}.  
\begin{lemma}\label{l2.2}
Notations and hypotheses as above. Then, the extension $k(X)\mid k(y)$ is a
Galois cyclic extension of degree $m$.
\end{lemma}
\begin{proof} From the proof of Proposition 2.1 we see that the extension
$k(X)\mid k(y)$ is ramified exactly at the rational points of Type 2 and that
$T_2 = (q+1)$. Moreover, this extension is totally ramified at those points.
Viewing the function $y$ on $X$ as a morphism of degree $m$
\[
y\colon X \to \mathbb P^1(\overline k),
\]
we consider the elements of the finite field $k$ belonging to the set below
\[
V_1 = \{y(P) \mid P \in X(k) \quad \text{is of Type 1}\}.
\]
Since $T_2=(q+1)$  and the rational points of Type 2 are totally ramified for
the morphism $y$, we must have $\#\,V_1 \le (q^2-q)$. Now, above each
element of the set $V_1$ there are at most $m$ 
rational points of the curve $X$, those points being necessarily of Type 1,
and hence:
\[
\#\,X(k) = q^2+1+2gq = T_1+T_2 \le m\cdot (q^2-q)+(q+1).
\]
From the genus formula in Proposition 2.1, we then conclude that $\#\,V_1 =
(q^2-q)$ and also that above each element of the set $V_1$ there are exactly
$m$ rational points of the curve.

Now the proof continues as in the proof of \cite[Lemma 5]{R-Sti}. We will 
repeat their argument here for completeness. Let $\widetilde F$ be the 
Galois closure of the extension $k(X)\mid k(y)$. The field $k$ is still
algebraically closed in $\widetilde F$ since the elements of the set $V_1$
split completely in $k(X)\mid k(y)$. Moreover the extension $\widetilde
F\mid k(X)$ is unramified, as follows from Abhyankar's lemma
\cite[ch.III.8]{Sti}. Hence,
\[
2\tilde g-2 = [\widetilde F\colon F](2g-2),
\]
where $\tilde g$ denotes the genus of the field $\widetilde F$. The 
$(q^2-q)$ elements of the set $V_1$ split completely in $\widetilde F$ and 
then they give rise to $(q^2-q)m [\widetilde F\colon F]$ rational points 
of $\widetilde F$ over $k$. Then, from the Hasse-Weil bound, we
conclude
\[
(q^2-q)m[\widetilde F\colon F] \le q^2+2q+(2\tilde g-2)q =
q^2+2q+[\widetilde
F\colon F](2g-2)q.
\]
Substituting $2g = (m-1)(q-1)$ in the inequality above, we finally get:
\[
[\widetilde F\colon F] \le \frac{q+2}{q+1} \text{ and hence } \widetilde
F=F.
\]
Note that the extension is cyclic since there exist rational points (those 
of Type 2) that are totally ramified for the morphism $y$.
\end{proof}
\begin{theorem}\label{t2.3}
Let $X$ be a maximal curve of genus $g$ such that there exists a rational point
$P_0 \in X(k)$ with $m\cdot n = (q+1)$, where $m$ is a non-gap at $P_0$\,. Then
the curve $X$ is $k$-isomorphic to the curve given by the equation:
\[
y_1^q + y_1 = x_1^m\,.
\]
\end{theorem}

\begin{proof}
We know that $k(X)\mid k(y)$ is a Galois cyclic extension of degree $m$ and
moreover that the functions $1,y,y^2,\ldots,y^{n-1}$ and $x$ form a basis for
$L(qP_0)$. Let $\sigma$ be a generator of the Galois group of $k(X)\mid k(y)$.
Since $P_0$ is totally ramified, then $\sigma(P_0)=P_0$ and hence
$\sigma(L(qP_0)) = L(qP_0)$. Note that the functions $1,y,y^2,\ldots,y^{n-1}$
form a basis for the subspace $L((n-1)m\,P_0)$ and that $\sigma$ acts as the
identity on this subspace. Since $m$ and $q$ are relatively prime, we can
diagonalize $\sigma$ on $L(q\,P_0)$. Take then a function $v \in L(q\,P_0)$, $v
\notin L((n-1)m\,P_0)$, satisfying $\sigma(v) = \lambda v$, with $\lambda$ a
primitive $m$-th root of 1. 

\noindent Then denoting by $N$ the norm of $k(X)\mid
k(y)$, we get
$N(v) = (-1)^{m+1}\cdot v^m$.

\noindent Hence $v^m \in k(y)$ and since it has poles only at $P_0$\,, we must
have $v^m = f(y) \in k[y]$. Since $\text{div}_\infty(v) = \text{div}_\infty(x)
=q\,P_0$\,, we see that $\text{deg }\,f(y)=q$. Now from the fact that there are
exactly $(q+1)$ totally ramified points of $k(X) \mid k(y)$ and that all of
them are rational, we conclude that $f(y) \in k[y]$ is separable and has all
its roots in $k$. After a $k$-rational change of coordinates, we may assume
that $f(0)=0$. Then, we get the following description for the set $V_1$\,:
$V_1 = \{\alpha \in k \mid f(\alpha) \ne 0\}$. Knowing that all points of
$X$ above $V_1$ are rational points over
$k$ and from the equation $v^m = f(y)$, we get:
\[
f^n(y) \equiv f^{nq}(y)\,\,\,\,\text{mod }(y^{q^2}-y). \qquad\qquad\qquad(*)
\]

\noindent{\bf Claim}. $f(y) = a_1\,y+a_q\,y^q$, with $a_1,a_q \in k^*$.

We set $f(y) = \sum\limits_{i=1}^q a_i\,y^i$
and $f^n(y) = \sum\limits_{i=n}^{nq} b_i\,y^i$. Clearly, the fact that $a_1,a_q
\in k^*$ follows from the fact that $f(y)$ is separable of degree $q$. Suppose
that the set $I$ below is non-empty
\[
I = \{2 \le i \le q-1 \mid a_i \ne 0\},
\]
and then define
\[
t = \min\,I \quad\text{and}\quad j = \max\,I.
\]
Clearly, we have 
$b_{(n-1)q+j} = n\cdot a_q^{n-1}\cdot a_j \ne 0$. Since the unique
solution for $i$ in the congruence\,\,\, $i\,q
\equiv (n-1)q+j
\text{ mod }(q^2-1)$ , \,\,\,$i$ being smaller than $q^2$, is the one given by
$i = (n-1) + j\,q$, it follows from ($*$) above that $b_{(n-1)q+j} =
b_{(n-1)+j\,q}^q \ne
0$.

\noindent It now follows that $\text{deg }(f^n(y)) = n\,q \ge (n-1)+j\,q$
and hence we get that $n-j \ge 1$ if $n \ge 2$. Note that if $n=1$, then we get
$j \le 1$ and the proof of the claim is complete in this case. From now on we
then assume $n\ge2$. We then conclude that $t \le j \le (n-1)$. Note that then
$(n-1)+t < q$. 

Clearly, we also have 
$b_{n-1+t} = n\cdot q_1^{n-1}, a_t\ne 0$. Since the unique solution for 
$i$ in the congruence\,\,\,$i\,q
\equiv n-1+t$\,\,$\text{mod }(q^2-1)$, \,\,\, $i$ being smaller than
$q^2$, is 
the one given by $i = (n-1+t)q$, it follows from ($*$) above that
$b_{n-1+t} = b_{(n-1+t)q}^q \ne 0$.

\noindent As before, it now follows that $n\,q \ge (n-1+t)q$, and hence $t \le
1$. This gives the desired contradiction and hence the set $I$ is empty,
thereby proving the claim.

Now we are in a position to finish the proof of the theorem. Denoting
\[
f(k) = \{f(\alpha) \mid \alpha \in k\} \quad\text{and}\quad H = \{\beta^m \mid
\beta \in k\},
\]
we have that $\mathbb F_q^*$ is a subgroup of $H\backslash\{0\}$ of index equal
to $n$. Moreover, using the fact that above $V_1$ there are only rational
points, we have:
\[
f(k) \subseteq H = \underset{\ell=0}{\overset{n-
1}{\bigcup}}\,\xi^{\ell\cdot m}\,\mathbb F_q\,,
\]
where $\xi$ denotes a primitive element of the field $k$; i.e., $\xi$ is a
generator for the multiplicative cyclic group $k^*$. Since $f(k)$ is a $\mathbb
F_q$-linear subspace of $k$ as follows from the above claim, we conclude that
its dimension is one and hence
that $f(k) =\xi^{r\cdot m}\,\mathbb F_q$\,, for some $r$. Finally, putting $x_1
= \xi^{-r}\,v$ and $y_1 = \epsilon\,y$, where $\epsilon$ is the unique element
of $k^*$ satisfying
\[
\text{Tr }(\epsilon\,\alpha) = \xi^{-r\cdot m}\,f(\alpha), \quad \forall\,
\alpha \in k,
\]
we conclude the proof of the theorem (Tr being the trace operator in 
$k(X) \mid k(y)$).
\end{proof}
\begin{remark}
Notations being as above. Suppose that $m\cdot n \le q+1$ ($m$ being a non-gap
at some rational point $P_0$ of $X$). Then, we have
$q+1 \ge m\cdot n \ge m_n(P_0) = q$, where the last equality follows from
Proposition 1.5(iv). In case
that $m\cdot n=q$, we conjecture that $2g = (m-1)q$ and the curve is 
$\mathbb
F_{q^2}$-isomorphic to a curve given by
\[
F(y) = x^{q+1},
\]
where $F(y)$ is a $\mathbb F_p$-linear polynomial of degree $m$. We have
not been able to prove this possible result yet.

We notice that if one could show that the morphism $\pi: X\to \mathbb
P^{n+1}$ is a closed embedding, then by Proposition \ref{p1.9} we would
have the claimed formula for $g$.

Finally we also notice that $(m_1(P)-1)q/2$, $P\in X(k)$, is an upper
bound for the genus of maximal curves. This follows from \cite[Thm. 
1(b)]{Le}.
\end{remark}
\begin{example}\label{e2.4}
There exist maximal curves that do not satisfy 
the hypothesis of Theorem \ref{t2.3}. We give two such examples below:

(i) Let $X$ be the maximal curve over $\mathbb F_{25}$ and genus $g=3$
listed by Serre in \cite[\S4]{Se}. Let  $m,5,6$ 
be the first three non-gaps at $P\in X(\mathbb F_{25})$. Here we have
$6P_0=g^3_6$. We claim that
$m=4$ (and so $nm>q+1$). Indeed, if $m=3$ by Proposition 2.1 
we would have $g=4$. 

This example also shows a maximal curve where all the rational points are
non-Weierstrass points: in fact, since $5 = {\rm char}(k)>2g-2$ the curve
is classical.

(ii) Let $X$ be a maximal curve over $\mathbb F_{q^2}$ of genus $g$.
Suppose that $q\ge 2g+2$ (e.g. the maximal curves in Proposition 
\ref{p?} here, \cite[Thm. 3.12, Thm.
3.16]{G-Vl}, \cite[Thm. 1]{I}). Then $X$ does not satisfy the hypothesis
of Theorem 2.3. In fact,
for
$P_0\in X(k)$ we have $m_{g+i}(P_0)=2g+i$ and then
$n=q-g$. Therefore $m_1(P_0)n\ge 2n\ge q+2$, the last inequality following from
$q \ge 2g+2$.
\end{example}
\section{Maximal curves of genus $(q-1)^2/4$.} 
As an interesting application of the preceding section we prove:
\begin{theorem}\label{t3.1}
Let $X$ be a maximal curve over $\mathbb F_{q^2}$ of genus $g = (q-1)^2/4$\,.
Then the curve $X$ is $\mathbb F_{q^2}$-isomorphic to the one given by
\[
y^q + y = x^{q+1/2}\,.
\]
\end{theorem}

\begin{proof}
From Equation (\ref{eq1.4}) which is Castelnuovo's genus bound applied to
the
linear system
$|(q+1)P_0|$, we have $n=2$ (see remark before Proposition 1.5). From Theorem
2.3, it suffices to prove the existence of a rational point $P$ over $k$ with
$m_1(P) = (q+1)/2$. This is clearly true for $q=3$ (since $g=1$ in this case)
and hence we can assume $q \ge 5$.

We prove firstly some lemmas:
\begin{lemma}\label{l3.2}
Let $P$ be a rational point over $k$ of the curve $X$ (hypothesis being as in
Theorem 3.1). Then we have that\,\, $\ell(2(q+1)P) = 9$ and that either $m_1(P)
=(q-1)$ or $m_1(P) = (q+1)/2$. Moreover, the divisor $(2g-2)P$ is a canonical
divisor.
\end{lemma}
\begin{proof}
Let $m_i = m_i(P)$ be the $i$-th non-gap at the rational point $P$. We have the
following list of non-gaps at $P$:
\[
0 < m_1 < m_2 = q < m_3 = q+1 \le 2m_1 < m_1+m_2 < m_1+m_3 \le 2m_2 < m_2+m_3 <
2m_3\,.
\]
The inequality $2m_1 \ge (q+1)$ follows from the fact that $n=2$ and $q$ 
odd. Clearly,
\begin{align*}
&m_3 = 2m_1 \quad\text{if and only if}\quad m_1 = (q+1)/2; \quad\text{and}\\
&m_1+m_3 = 2m_2 \quad\text{if and only if}\quad m_1 = (q-1).
\end{align*}
Since $q \ge 5$, one cannot have both equalities above simultaneously. From the
above list of non-gaps at $P$, it then follows that $\ell((2q+2)P) \ge 9$.
Moreover, after showing that $\ell((2q+2)P) = 9$, it also follows that either
$m_1 = (q+1)/2$ or $m_1 = (q-1)$. Let $\pi_2\colon X \to \mathbb P^{r+1}$ be
the morphism associated to the linear system $|(2q+2)P|$; we already know that
$r \ge 7$ and we have to show that $r=7$. Castelnuovo's bound for the morphism
$\pi_2$ gives
\[
2g = \frac{(q-1)^2}{2} \le M\cdot(d-1-(r-e)),\qquad\qquad (*)
\]
where $d = 2q+2$, $M = \bigg[\dfrac{d-1}{r}\bigg]$ and $d-1 = M\cdot r+e$.
Since $(r-e) \ge 1$ we have $d-1-(r-e) \le 2q$, and hence
\[
(q-1)^2 \le 4qM \quad\text{and then}\quad q^2-q \le 4qM,
\]
since the right hand side above is a multiple of $q$. For $r \ge 9$, we now see
that
\[
q-1 \le 4M \le 4\cdot \frac{2q+1}{9}\,, \quad\text{and then}\quad q \le 13.
\]
The cases $q \le 13$ are discarded by direct computations in Equation
$(*)$
above, and hence we have $r \le 8$. Now we use again Equation $(*)$ to
discard
also the possibility $r=8$. Since $q$ is odd, we have

\begin{align*}
2q+1 \equiv 3(\text{mod}\,8) \quad &\text{or}\quad 2q+1 \equiv
7(\text{mod}\,8).\\
\intertext{It then follows}
\begin{cases}
M = (q-1)/4\\
\text{and } e=3
\end{cases} \qquad &\text{or}\qquad
\begin{cases}
M = (q-3)/4\\
\text{and } e=7.
\end{cases}
\end{align*}
Substituting these two possibilities above in Equation $(*)$, one finally
gets
the desired contradiction; i.e., one gets
\[
(q-1)^2 \le (q-1)(q-2) \quad \text{or}\quad (q-1)^2 \le (q-3)\cdot q.\qquad
\]

Now we prove the last assertion of the lemma. One can easily check that both
semigroups $H_1$ and $H_2$ below are symmetric, with exactly $g=(q-1)^2/4$
gaps:
\[
H_1 = \langle(q-1),q,q+1\rangle \quad\text{and}\quad H_2 = \big\langle
\frac{q+1}{2},q\big\rangle.
\]
At a rational point $P$ on $X$ the Weierstrass semigroup $H(P)$ must then be
equal to $H_1$ or $H_2$\,. Hence the semigroup $H(P)$ is necessarily symmetric
and the last assertion follows.
\end{proof}

\begin{lemma}\label{l3.3}
Let $\d = |(q+1)P|$ with $P$ being a rational point of $X$ (hypothesis as in
Theorem \ref{t3.1}). Then at any non-rational point $Q$ of $X$, the 
$(\d,Q)$-orders are $0,1,2,q$. In particular the order sequence for $\d$ is
$0,1,2,q$, and the set of $\d$-Weierstrass points is exactly the set of
rational points.
\end{lemma}

\begin{proof}
Let $0,1,j,q$ be the $(\d,Q)$-orders. Consider the following set $S$:
\[
S = \{0,1,2,j,j+1,2j,q,q+1,q+j,2q\}.
\]
The set $S$ consists of $(2\d,Q)$-orders, and hence from Lemma \ref{l3.2} we
must have $\#\,S \le 9$. This eliminates the possibilities
\[
3 \le j \le (q-1)/2 \quad\text{and}\quad \frac{q+3}{2} \le j \le q-2,
\]
and it then follows that $j \in \{2,(q+1)/2,q-1\}$. From Lemma \ref{l3.2} we
know that
\[
(2g-2)P = \frac{(q+1)(q-3)}{2} P \quad\text{is canonical.}
\]
Then the following set $S(j)$ consists of orders at $Q$ for the canonical
morphism
\[
S(j) = \{a+bj+cq \mid a,b,c\in\mathbb N \text{ with } a+b+c \le
\frac{q-3}{2}\}.
\]
One can check that $\#\,S(j) = (q-1)^2/4$ if the value of $j$ belongs to
$\{2,(q+1)/2, q-1\}$, and hence that $S(j)$ consists of all canonical orders at
the point $Q$. Then the set $H(j)$ below is necessarily a semigroup:
\[
H(j) = \mathbb N \setminus (1+S(j)).
\]
This semigroup property on $H(j)$ is only satisfied for the value $j=2$, as one
checks quite easily, and this finishes the proof of this lemma.
\end{proof}

Now we turn back to the proof of Theorem 3.1. Suppose that $m_1(P) = (q-1)$ at
all rational points $P$ on the curve. It then follows from Proposition 1.5(iii)
that the $(\d,P)$-orders are $0,1,2,q+1$ and hence $v_P(R) = 1$, where $R$ is
the divisor supporting the $\d$-Weierstrass points. On the other hand, we have
\[
\text{deg }R - \#\,X(k) = 3(2g-2) - (q-3)(q+1) = \frac 12 (q+1)(q-3).
\]
Since $q \ge 5$ and $v_P(R)=1$ for $P$ rational, we would then conclude the
existence of non-rational points that are $\d$-Weierstrass points. This
contradicts Lemma \ref{l3.3} and hence, from Lemma \ref{l3.2}, we finally
conclude the existence of a rational point $P$ satisfying
\[
m_1(P) = (q+1)/2. 
\]
\end{proof}
We can explore further the idea of the above proofs to obtain a partial
analogue of the main result of \cite{F-T}, namely
\begin{scholium}
Let $X$ be a maximal curve over $k$ whose genus $g$ satisfies
$$
(q^2-3q+2)/4<g\le (q-1)^2/4.
$$
If $q$ is odd, neither $q$ is a power of 3 nor $q\not\equiv 3\pmod{4}$,
then $g=(q-1)^2/4$.
\end{scholium}
Notice that Example 2.4(i) shows that the hypothesis on $g$ above is
sharp. This Scholium is the first step toward a characterization of a
maximal curve whose genus is $\frac{q-1}{2}(\frac{q+1}{t}-1)$ with $t\ge
3$.

\end{document}